\documentclass{article}
\usepackage[utf8]{inputenc}
\usepackage[a4paper, total={6.15in, 9.7in}]{geometry}
\usepackage{amsmath}
\usepackage{mathtools}
\usepackage{amsthm}
\usepackage{amssymb}
\usepackage[english]{babel}
\usepackage{dirtytalk}
\usepackage{graphicx}
\usepackage{authblk}
\usepackage[hyphens]{url}
\usepackage{hyperref}
\usepackage{caption}
\usepackage{subcaption}
\usepackage{cases}
\usepackage{xcolor} 
\usepackage{bm}
\usepackage{tikz}
\usetikzlibrary{shapes,arrows.meta,positioning}

\usepackage{tabularx}
\usepackage{blkarray}
\usepackage{enumerate}
\usepackage{comment}
\usepackage{algorithm}
\usepackage{algorithmic}
\usepackage{float}
\usepackage{booktabs}
\usepackage{ragged2e}

\usepackage{appendix}
\usepackage{multirow}

\title{Optimal Control of a Battery Storage On the Energy Market}
\author[1]{Stephan Schl{\"u}ter}
\author[2]{Abhinav Das}
\author[3]{Matthew Davison}
\affil[1]{Institute of Energy Engineering and Energy Economics, Ulm University of Applied Sciences}
\affil[2]{Faculty of Mathematics and Economics, Ulm University}
\affil[3]{Departments of Mathematics and of Statistical \& Actuarial Sciences, Western University London, ON, Canada}

\date{}

\begin{document}
\maketitle

\begin{abstract}
Electricity storage is crucial for a successful transition towards carbon-neutral energy production. Despite considerable research and a number of promising future alternatives such as hydrogen, battery storages currently remain the first choice. However,  costs remain high and it remains to be shown whether an investment can be profitable. This article addresses this question by modelling a battery storage operating in the German power market. We consider two periods with very distinct price dynamics, namely a calm year (2020) and a turbulent year (2023). It shows that even for low battery costs a 2020 style price environment does not allow for profitable battery operation, whereas current market conditions allow for positive payoffs.
\end{abstract}

\section{Introduction}

In Germany in 2024, renewable energy is available in considerable quantity and at reasonable prices. The stochastic supply dynamics of renewable power means that conventional power plants must sometimes ramp up to meet relatively inelastic demand, inducing additional costs and unwanted carbon dioxide emissions. The solution to this dilemma is to store power, and there has been substantial recent research in this field. For example, Sharma and Mortazavi \cite{sharma2023pumped} and  Geth et al. \cite{GETH20151212} discuss how to overcome geographical constraints for pumped thermal electricity storage and pumped hydro storage. In practice, there are currently three main options: heat, pumped hydro storage  and electric batteries. Heat is commonly used at small scale in households while the others have been installed at larger scale \cite{website1,GETH20151212}. Pumped hydro energy storage is widely used worldwide although, unlike batteries, it  faces geographical constraints. This article focuses on exploring the profitability of battery storages which,  despite recent price drops, remains expensive. An active operation strategy is needed to ensure the profitability required to ensure that sufficient battery storage is built for the cost- and emissions- efficient use of renewable energy. We focus on deriving an optimal strategy for operating a battery storage using Germany as a sample power market. Based on our findings, in a next step, renewable sources can also be connected to the battery.

A battery can be considered as a real option, the optimal operation of which is the subject of a rich literature. Boogert \& de Jong use the least squares Monte Carlo (LSMC) method \cite{boogert2006gas} to price a gas storage while Thompson et al. \cite{thompson2009natural} propose a  partial differential equations-based method, with recent efforts to use machine learning methods as well \cite{wu2023multi}. Batteries differ both in cost and in physical properties from other energy storage modes. They have a much higher cycling frequency (seasonal battery storages are rare) and have high charging and discharging rates, so it is reasonable to assume that a full charge/discharge cycle is possible within a single day. This motivates us to focus on a battery optimization horizon of a single week, within which hourly day ahead spot prices are considered. Thereby we find that complexity can be significantly reduced and propose a new efficient model by discriminating between intra-  and inter-day optimization. Finding that the profit of intraday optimization closely resembles a linear function, we simplify the inter-day optimization: Based on forecasts generated by a vector autoregressive (VAR) model we determine the optimal storage level for the end of the upcoming day, which in turn allows an intraday optimum to be found. We test our strategy on two years with very distinct price behaviour: the low price, low volatility Year 2020 and the high price, high volatility Year 2023. In this way we see that the battery's profitability is higher with volatile prices, consistent with economic theory. Profitability also increases with battery size and charge/discharge rates, as these allow these volatile prices to be financially exploited. Third, the comparably simple VAR model is sufficient as intraday optimization is by far more important than inter-day decisions, i.e. storage levels at the end of the day. In summary, we prose a new lean and efficient model for the valuation of a battery storage. Contrary to current literature (\cite{feng2022optimization, zhang2021arbitrage}), which concerns with handling a 24 dimensional price vector we are able to reduce the model's dimension towards deciding about an univariate variable, namely the storage level at the end of the day. Moreover, we extend the discussion of storage profitability as we not only analyze cash flows but also consider profitability in terms of time to break even.

The article is structured as follows: In Section \ref{Sec:priceForecast} we comment on price forecasting, in Section \ref{sec:electricityControl} we present or optimal control approach. Section \ref{sec:caseStudy} contains the case study and Section \ref{sec:conclusion} concludes the paper.

\section{A Power Price Forecasting Model}
\label{Sec:priceForecast}
The scope of operation for the battery here is one week, and we consider hourly day ahead prices. In the German power market, Saturday, Sunday, and Monday are simultaneously traded on Friday.  This feature has been ignored in almost all existing price models;  we will also ignore it in our price forecasts. However, below, we will consider this fact when developing our control strategy. For simplicity, we use a vector autoregressive model of order one, VAR(1), which sets the 24-dimensional hourly price vector in linear relation to the previous day's prices incorporating an error vector. Calibration is done using maximum likelihood estimation \cite{lutkepohl2005new}. For a more detail explanation please see  Appendix \ref{appendix_VAR}. Let $\vec{P}^t = P_i^t, i = 1,\ldots,24$ denote the vector of power prices with delivery on day $t$ and assume that the error vector $\vec{\epsilon}^t$ is Gaussian distributed with zero mean and $24 \times 24$ dimensional covariance matrix $\Sigma$. Then the VAR(1) model reads as follows:
\begin{equation}
\vec{P}^t = A \vec{P}^{t-1} + \vec{\epsilon}^{\hspace{0.7mm}t},
\label{eq:VAR}
\end{equation}
with $A \in \mathbb{R}^{24 \times 24}$. For details on model calibration and application, please refer to \cite{embrechts2011quantitative}. Although the VAR model offers only a mediocre forecasting quality, it is comparably easy to calibrate and to handle. G\"{u}lerce \& Gazanfer \cite{MR3786477} generalize the model by including a moving average part and wavelet coherence. Other authors allow for stochastic volatility \cite{MR4667861} or apply advanced machine learning methods \cite{MR4087749}.

\section{Optimal Battery Storage Control}
\label{sec:electricityControl}

There is ample research and literature regarding the optimal control of an energy storage \cite{wu2023multi} \cite{aazami2022optimal}, \cite{anderson2015optimal}. Articles vary in scope. Some focus on the physical side, some concentrate on the randomness of the input such as wind. E.g., Hannan et al. \cite{HANNAN2021103023} analyze battery system optimization, addressing objectives such optimal sizing, and systematically categorizing them. Babatunde et al. \cite{https://doi.org/10.1002/er.4388} integrate demand side management with the storage optimization. Castillo and Gayme \cite{CASTILLO2014885} explore the potential and obstacles of grid-scale energy storage while considering various storage technologies.  Anderson Burke \& Davison optimize a stylized storage in face of random and intermittent production \cite{ABD2015} Finally, a comprehensive review of a wide range of articles is provided by \cite{machlev2020review}. In this article we focus on the financial benefit of using a battery storage in developing market trading strategies. 

In Germany hourly prices are fixed day ahead at noon, hence the intraday price structure of tomorrow is known and finding the maximum profit means solving a linear optimization problem subject to the storage level at both the beginning and end of the day (see Section \ref{SecIntraday}). The corresponding inter-day optimization is described in Section \ref{sec:interday}.

\subsection{Intraday Optimization}
\label{SecIntraday}

Let $\vec{x} = x_i, i = 1,\ldots,24$ denote the charging and discharging amounts per hour. We adopt the convention in which a negative value means discharging and a positive value means charging the battery. Let $x_{max}^{out} < 0 $ be the maximum discharge amount per hour and $x_{max}^{in}$ the maximum charging speed per hour, respectively. Let $B^{max}$ be the maximum storage charging level and $0 \leq B^t \leq B^{max}$ the charging level at the end of Day $t$. Then the linear optimization problem reads as follows:
\begin{equation}
\label{eq:intraday}
\max_{x_i}{\sum_{i = 1}^{24}{-x_i P_i}}.
\end{equation}
subject to:
\begin{eqnarray*}
\sum_{i=1}^{24}{x_i} & = & B^t - B^{t-1} \\
x_{max}^{out} \leq & x_i  & \leq x_{max}^{in} \\
B^{t-1} +  \sum_{i=1}^{j}{x_i} & \geq & 0, \; \; \; \hbox{for} \; j = 1,2,\ldots,24. \\
B^{t-1} +  \sum_{i=1}^{j}{x_i} & \leq & B^{max}, \; \; \; \hbox{for} \; j = 1,2,\ldots,24. \\
\end{eqnarray*}
Note that we have an additional feasibility constraint for $B^{t}$ and $B^{t-1}$, i.e. $B^{t}-B^{t-1}$ has to be possible given the storage fundamentals:
\begin{equation*}
24 x_{max}^{out} \leq B^t - B^{t-1} \leq 24 x_{max}^{in}.
\end{equation*}
\subsection{Inter-Day Optimization}
\label{sec:interday}
\begin{figure}[htbp]
\centerline{\includegraphics[scale=0.4]{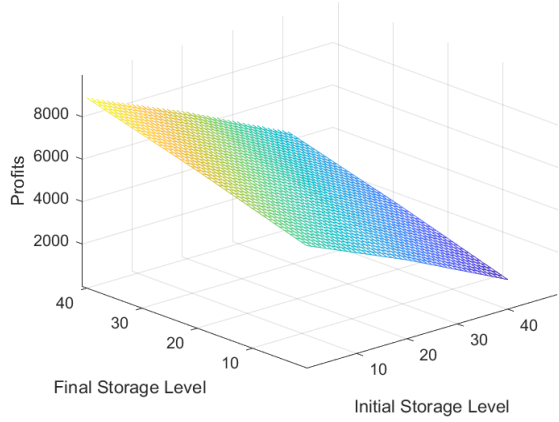}}
\caption{Values for the function $Profits$ depending on different values for the initial and final storage levels $B^{t-1}$ and $B^t$}.
\label{planarProperty}
\end{figure}

Let the optimal solution to Eq. \ref{eq:intraday} be denoted by the function $Profits(B^{t-1},B^t,\vec{P}^t)$. Note that now the inter-day optimization is reduced to deciding the storage level at the end of Day $t$.  Given a price vector $\vec{P}^t$ this function is nearly a linear plane, considerably simplifying the optimization problem: We can reasonably approximate the optimization function as shown in Eq. \eqref{eq:optim}.
\begin{eqnarray}
\label{eq:optim}
\max_{B_t}{E \left(\sum_{i = t}^T{Profits(B^{i-1},B^{i},\vec{P}^t)} \right)} \\ \nonumber
\approx \max_{B_t}{ \sum_{i = t}^T{ E \left(Profits(B^{i-1},B^{i},\vec{P}^t) \right)}}
\end{eqnarray}
where $T \in \mathbb{N}$ denotes the last day of the storage operation period and the expectation is considered to be over the multivariate price process $\vec{P}^t, \vec{P}^{t+1},\ldots,\vec{P}^T$. This means we need not consider the distributional properties of the price vector, a point forecast of the prices is sufficient. An example of the planar property of the function $Profits$ is shown in Fig. \ref{planarProperty}. Here the function values are shown for the hourly power prices of  $24^{\hbox{th}}$ May 2023 and all potential initial and final storage levels on this day. We see that the profit evolves almost linear in both inital and final storage level.

The optimal control at time instance $t$ is now determined in a two-step procedure. First, compute forecasts for the days $t+1,\ldots,T$ using one of the price models from Section \ref{Sec:priceForecast}. Second, based on these price forecasts, start at time $T$ with a predetermined $B^T = B^{max}$ and apply the principle of dynamic programming: Move backwards in time to determine the optimal storage level $B^t$ using the continuation function $C(\hat{\vec{P}}^{t+1},\ldots,\hat{\vec{P}}^{T} | B^{t})$, i.e. the expected payoff of tomorrow given a decision $B^t$ of today and price forecasts for time instances $t+1,\ldots,T$. In this way identify the optimal battery level for each time instance $t$ given an initial level $B^{t-1}$ of the respective day $t$ via Eq. \eqref{eq:bellmann}.
\begin{equation}
\label{eq:bellmann}
\max_{B_t}{\left[ Profits(B^{t-1},B^{t},\vec{P}^t) + C(\hat{\vec{P}}^{t+1},\ldots,\hat{\vec{P}}^{T}|B^{t})\right]}
\end{equation}

\section{Case Study}
\label{sec:caseStudy}
In this case study we consider the optimal control of an electricity storage with a capacity of $40$megawatt hours (MWh), about the size of a storage that has recently been built by the German power utility RWE \cite{website1}. We test two scenarios of battery power: In Scenario 1 the maximum charging speed $x_{max}^{in}$ is $ 20$MW and the maximum discharging speed of $x_{max}^{out}$ is $-20$MW, i.e. $50$\% of the battery capacity. In Scenario 2 we set $x_{max}^{in} = 5$MW and $x_{max}^{out} = -5$MW, i.e. $12.5$\% of the battery power. Prices are provided by the European Power Exchange (EPEX Spot) and our history ranges from January 2018 to Dec 2023 (Section \ref{dataSet}). We use the dataset to calibrate the VAR model. The model itself is tested only for the years 2020 and 2023 to compare the battery's profit in two different environments: one period of relatively calm prices (2020) and one volatile period (2023). So, for each year, we consider weekly blocks, i.e. we expect the battery to be full every Sunday midnight. Hence the battery needs to be optimized only with a horizon $T$ in Eq. \ref{eq:optim} of 7. 

\subsection{The Dataset}
\label{dataSet}

\begin{figure}[htbp]
\centerline{\includegraphics[height = 4.7cm, width = 9.5cm]{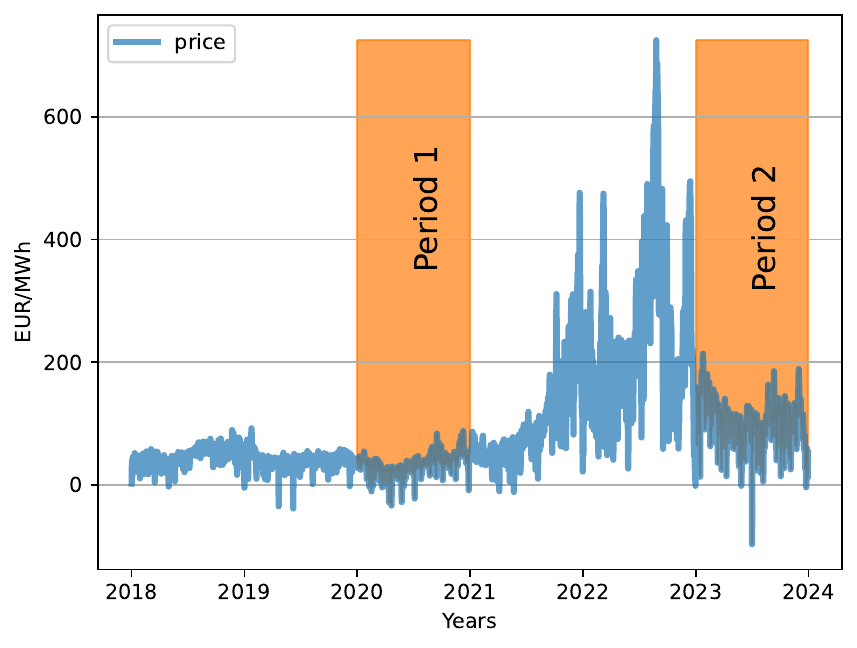}}
\caption{Daily Average over Hourly German Power Prices}
\label{fig:prices}
\end{figure}

In Fig. \ref{fig:prices} the daily average over the hourly German power prices are shown. In orange color, the two test periods, namely 2020 and 2023 are indicated. We can clearly see the gradual outbreak of the Ukraine crisis in summer 2021 caused by a considerable reduction of gas flows from Russia to Germany. Prices peaked in 2022 and the situation eased with increasing liquidity on the power market in 2023. Nevertheless, both test periods, namely 2020 and 2023, differ considerable with regards to essential properties as can be seen in Fig. \ref{fig:prices}. In 2020 negative prices (even on a daily average level) were common whereas in 2023 prices were mainly positive. Volatility also varies considerably across periods, with a standard deviation of the first differences of 30 EUR/MWh in 2023 compared to 11.54 EUR/MWh in 2020. The Skewness  and excess of the daily averages are both small for each period. Hence, at least for the first differences of daily average prices, a Gaussian distribution cannot be discarded.

\subsection{General Results}
\label{ResultsSection}

In Table \ref{tab:annualPayoff} the annual payoffs for two different battery power levels and the years 2020 and 2023 is shown. We can see that in general profits are much lower in 2020 than in 2023. For example, in 2020 a battery with a power of 5MW generated only about one third of the profits that would have been possible in 2023. This corresponds with the difference in price volatility and absolute prices (see Fig. \ref{fig:prices}), which is also in line with traditional option pricing insights, which demonstrate the volatility is a significant driver of option value.  
\begin{table}[htbp]
\caption{Annual Payoff per 1000 MWh in EUR Depending on Battery Power}
\begin{center}
\begin{tabular}{|l|c|c|}
\hline
\textbf{Year} & \textbf{Power 5 MW} & \textbf{Power 20 MW} \\ \hline
\textbf{2020}                      &     7.318    &   13.61   \\
\textbf{2023}                      &     22.05    &   42.23   \\ \hline
\end{tabular}
\label{tab:annualPayoff}
\end{center}
\end{table}
\vfill

\begin{figure}[htbp]
\centerline{\includegraphics[height = 6.7cm, width = 9cm]{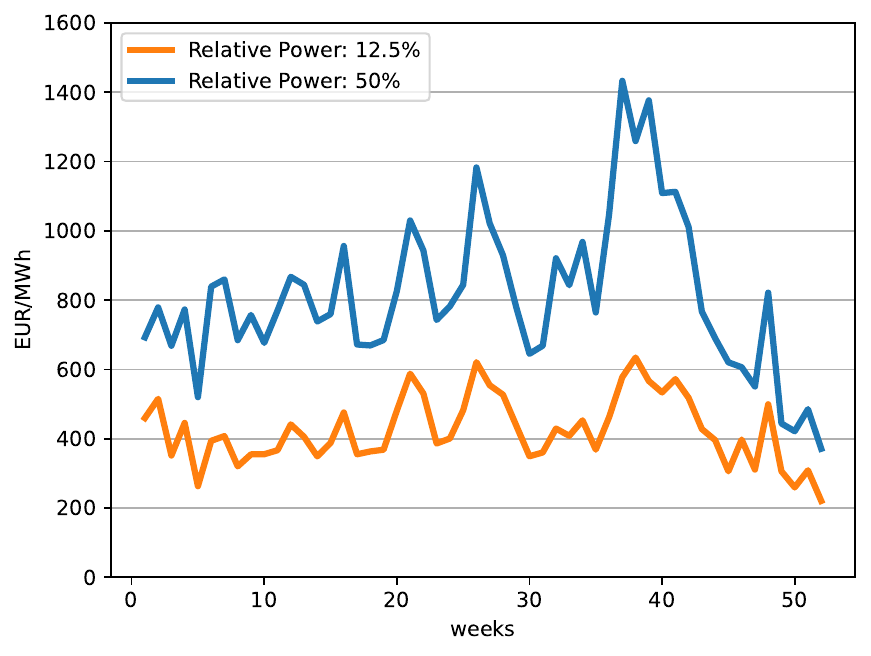}}
\caption{Weekly Payoff in 2023 for a 40MWh Battery and Two Different Power Levels} 
\label{fig:results2023}
\end{figure}

For 2023 we also depict the weekly payoff in Figure \ref{fig:results2023}, where the cumulative weekly payoff is shown for all weeks of 2023 and for two relative battery power levels. Thereby the relative power level refers to the battery's (dis)charging capability within one hour relative to its size. For example, the blue line indicates a comparably high relative power of $50$\% meaning  that the battery can be fully (dis)charged within two hours.In the graph we cannot see a distinct seasonal pattern but the influence of the battery power. The higher the power, the higher the intraday profit, which is reasonable: A battery with higher power can purchase more electricity at the cheapest hour and sell more at the most expensive hour. 

\begin{figure}[htbp]
\centerline{\includegraphics[height = 3cm, width = 7cm]{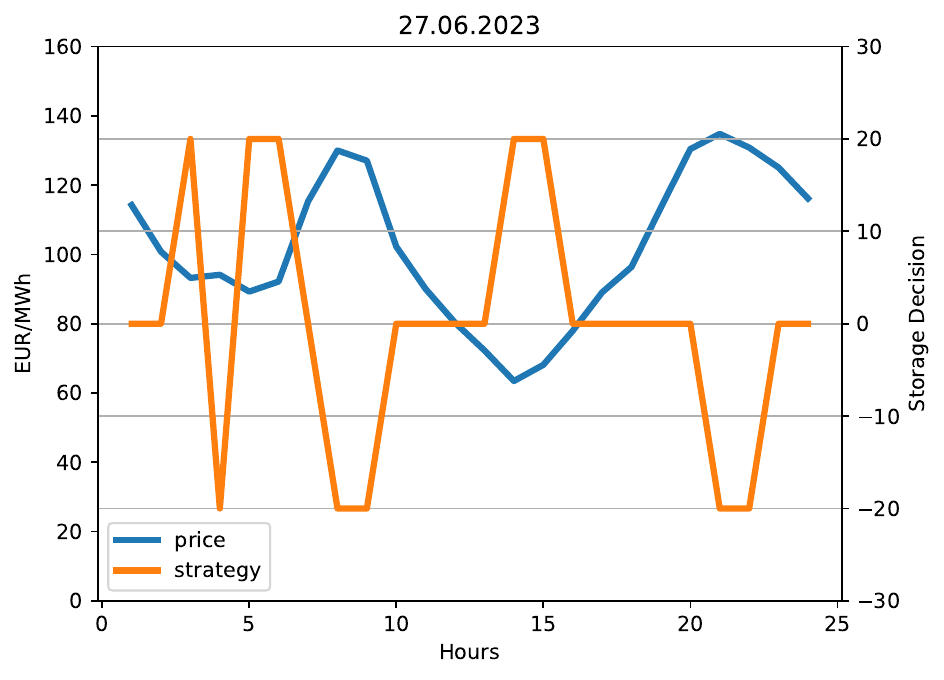}}
\centerline{\includegraphics[height = 3cm, width = 7cm]{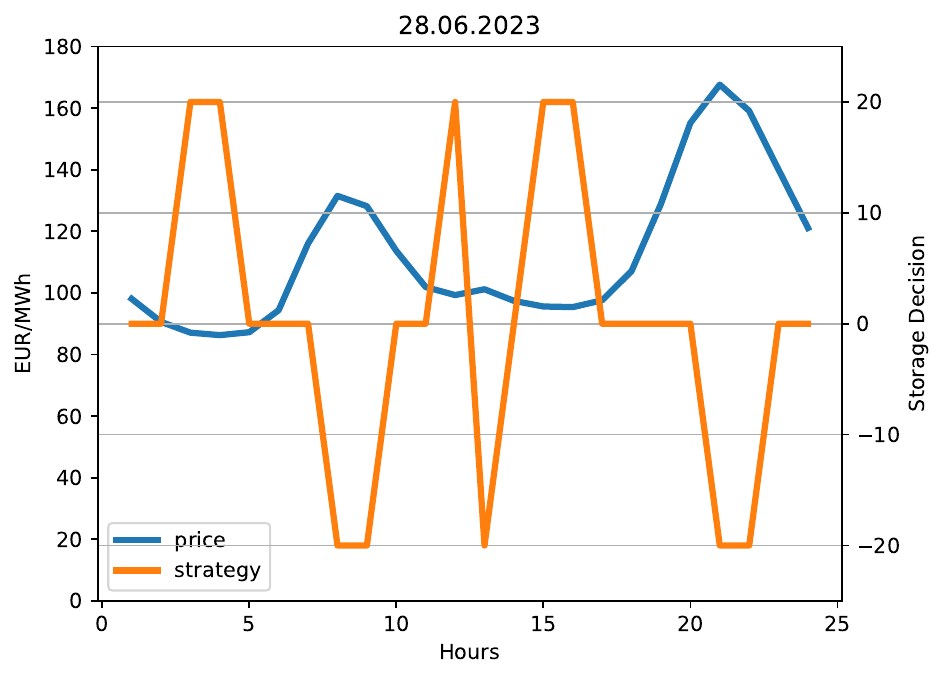}}
\caption{A Comparison of Hourly Prices vs. Optimal Stroage Strategy on Two Consecutive Days for 20MW Battery Power}
\label{fig:results2023_20}
\end{figure}

\begin{figure}[htbp]
\centerline{\includegraphics[height = 3cm, width = 7cm]{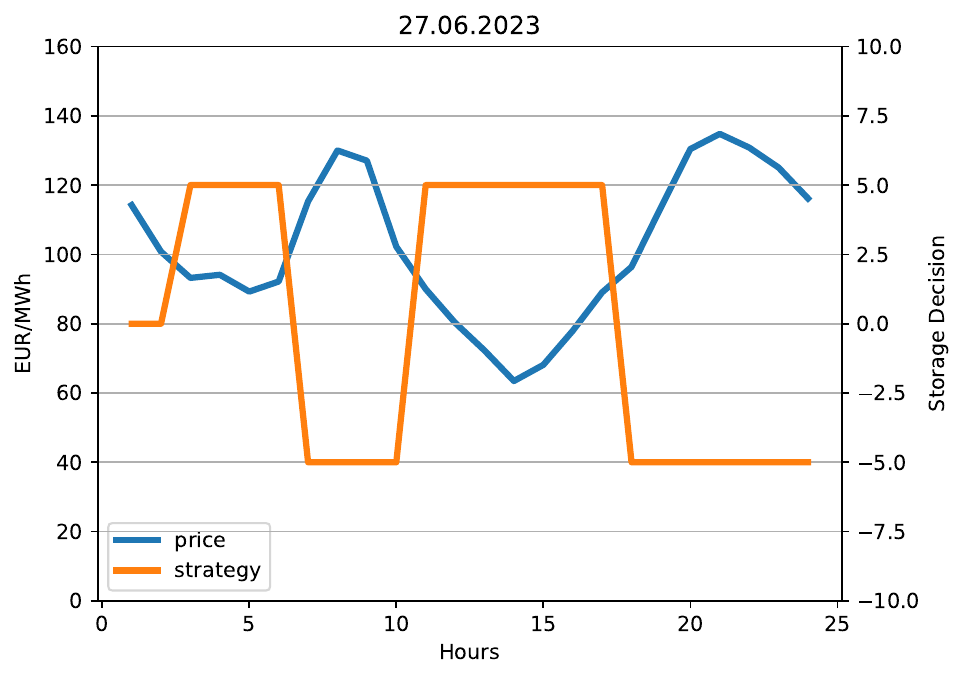}}
\centerline{\includegraphics[height = 3cm, width = 7cm]{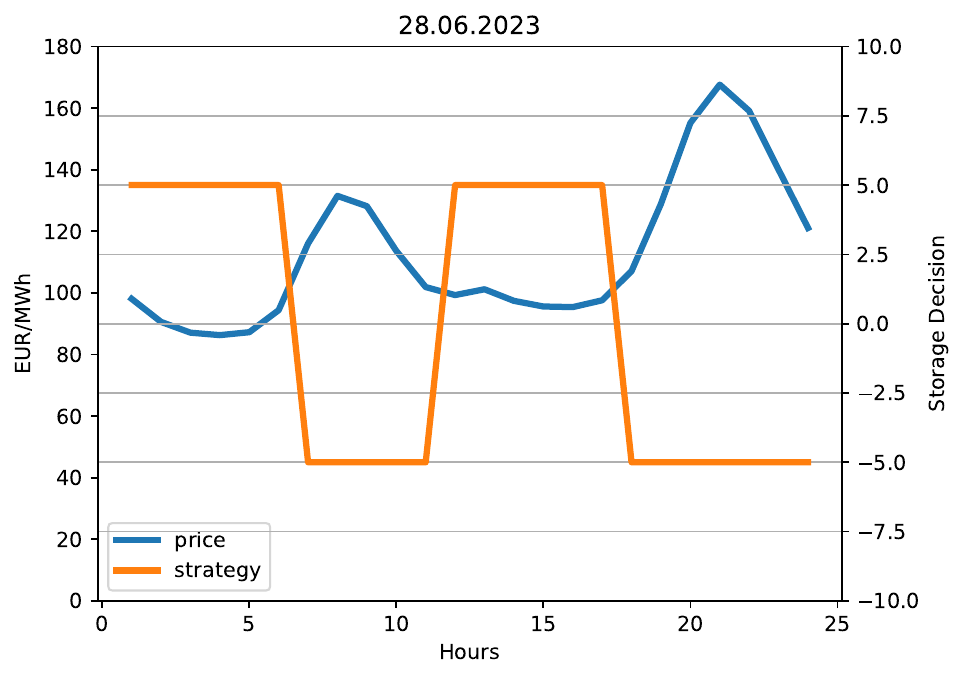}}
\caption{A Comparison of Hourly Prices vs. Optimal Storage Strategy on Two Consecutive Days for 5MW Battery Power}
\label{fig:results2023_5}
\end{figure}

Interestingly, for both years and both tested battery power, the results are close to the profits under perfect foresight, in which the true future prices are assumed to be known. This may seem remarkable as we apply only a comparably simple VAR model for forecasting the near future, however it shows that the intraday optimization dominates the inter-day optimization. This illustrates that no exact forecast is required, but only the rough price structure as intraday optimization means trading on price differences. Hence, we can conclude that a sophisticated forecasting model is not needed here. As an example we compared our algorithm's decision regarding the storage levels $B^t$ to the perfect foresight decision in Week 38 of 2023. Results are given in Table \ref{tab:decisionsB}, whereby numbers indicate the storage levels in Mwh at the end of the respective day.  Here, we differ only on two days from the optimal strategy.
\vfill
\begin{table}[htbp]
\caption{Optimal Storage Level (in MWh) at the End of a Trading Day}
\begin{center}
\begin{tabular}{|l|c|c|c|c|c|c|c|}
\hline
 & \textbf{Mon} & \textbf{Tue}  & \textbf{Wed}  & \textbf{Thu}  & \textbf{Fri}  & \textbf{Sat}  & \textbf{Sun} \\ \hline
\textbf{Our Strategy}    &    20    &  0    & 0    & 0    & 0    & 0    & 40         \\
\textbf{Perfect Foresight}  &    0    &  0    & 0    & 20    & 0    & 0    & 40       \\ \hline
\end{tabular}
\label{tab:decisionsB}
\end{center}
\end{table}
In order to shed more light on the control decision of our strategy, we exemplary show the results for two consecutive days in June 2023 in Figures \ref{fig:results2023_20} and \ref{fig:results2023_5} These graphics visualize the optimal strategy for charging and discharging the battery over a day via proposed optimization problem (\ref{eq:intraday}) for battery power 20 and 5 MW respectively. The blue line indicates the hourly price over the day, the orange line shows the trading decision. Positive values mean that the battery is charged, negative values indicate a discharging, i.e. selling decision. For example, on all days for both battery power versions, the algorithm charges the battery in the afternoon hours before selling at the highest possible price at the evening peak hours.


Eventually, to shed light on the influence of battery, power we test its impact on the weekly profits using Week 38 of 2023 as an example. In Figure \ref{fig:results2023} we see a sub-linear, potentially logarithmic effect, i..e relative power positively effects the profits, but with decreasing effect. The difference in profits between 5\% and 10\% is considerably larger than between 80\% and 85\%.

\begin{figure}[htbp]
\centerline{\includegraphics[height = 4.5cm, width = 8cm]{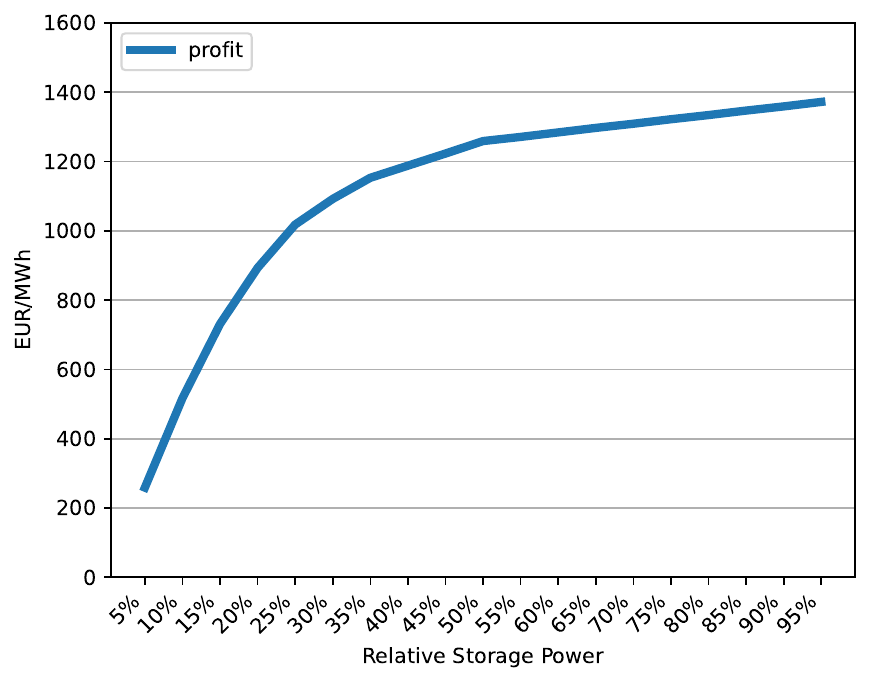}}
\caption{Relative Storage Power Vs. Profits per MWh for Week 38 in 2023}
\label{fig:results2023}
\end{figure}

\subsection{Investment Costs}
For evaluating a battery's profitability, both investment and operational costs are needed. However, recent market fluctuations make estimating these costs a challenge.  We follow Cole \& Karmakar \cite{cole2023cost} who expect  declining costs and compute the time to break even, i.e. the time needed until the investment object becomes profitable, for a range of costs. Results are summarized in Table \ref{tab:ttm}, where we show the time to break even in years for a range of battery cost levels, different battery power levels and the years 2020 and 2023 in order to consider two different market environments. In contrast to 2020 the results for 2023 are promising. In 2020, even in the best price scenario and for a high battery power, the time to break even amounts to 8 years, which -- depending on the properties -- might even approach the battery lifetime. All other financial scenarios clearly speak against battery storage economics -- especially those with 5 MW battery power. Here nearly all time to maturity numbers  are clearly beyond the battery lifetime. In the 2023 high power scenario a reasonable time to break even can be achieved even when investment costs are moderate. A 2021 report by the U.S. National Renewable Energy Laboratory \cite{cole2021cost} (NREL; Table 4 on Page 8) summarizes various recent battery cost estimates which show that prices of $300$EUR to $500$EUR remain more realistic.  However, this report projects that these prices will continue to decline rapidly to the $100$EUR to $250$EUR range by 2030, at which point break even times will be reasonable for the 2023 scenarios. 
\begin{table}[htbp]
\caption{Time to Break Even in Years Depending on Battery Power}
\begin{center}
\begin{tabular}{|l|c|c|c|c|c|}
\hline
\multicolumn{6}{|c|}{\textbf{Year 2020}} \\\hline
Power  & \multicolumn{5}{|c|}{Battery Costs per MWh (1000 EUR)} \\\hline
       &  100 & 200 & 300 & 400 & 500 \\ \hline
 5 MW  & 12   & 28   & 42  & 55 & 69  \\
 20 MW &  8  & 15   & 23 & 30  & 37 \\ \hline
 \multicolumn{6}{|c|}{\textbf{Year 2023}} \\\hline
Power  & \multicolumn{5}{|c|}{Battery Costs per MWh (1000 EUR)} \\\hline
       &  100 & 200 & 300 & 400 & 500 \\ \hline
 5 MW  &  5 &  10 &  14  &  18 & 23  \\
 20 MW & 3   &  5  & 8 &  10 & 12 \\ \hline
\end{tabular}
\label{tab:ttm}
\end{center}
\end{table}
\subsection{Policy Implications}
Batteries can help to buffer time-varying and partly unpredictable renewable power generation. 
Particularly from the perspective of European Union (EU), battery storage is a major issue. According to a report published by the European Union in 2021 \cite{european2020proposal}, EU is encouraging to scale up the battery production up to 19 times the today's battery production globally in order to reduce carbon emissions. This is only possible if both business users and private consumers are motivated to install batteries in considerable scale. As already mentioned in the introduction, costs are a major factor. This study puts light on how battery storage can be used efficiently and for energy arbitrage. Here, we provide a comparatively simple yet useful control strategy. Besides, it shows that given today's battery prices the goal of the EU is not realistic. Policymakers have to take measures accordingly, especially with a focus on reducing battery costs. Subsidies or tax exemptions would also reduce investment costs. 

\section{Conclusion}
\label{sec:conclusion}
In this article we operated a battery storage on the German power market trading hourly volumes. We tested two different periods (high volatility and low volatility) and two different battery powers, namely $50$\% of the capacity and $12.5$\% of the capacity. As a price model we used a simple VAR model. We also split the optimal control into intraday and inter-day optimization. In this way we saw that profits are comparably low in 2020, whereas the 2023 price level and volatility allows for a profitable storage operation -- depending on the investment costs.

However, this analysis does includes neither household consumption nor a residentially located renewable source such as wind or solar power, which is sometimes called prosumer model \cite{violi2023dealing}. Besides, a more detailed sensitivity analysis of the model's performance under different conditions is required. Future work needs to investigate how to integrate a battery on a in this context in a first step. In a second step,the scope has to be extended to larger grids such as streets or residential blocks with the purpose of stabilized the grid in presence of stochastic renewable energy by optimized battery control. However, combining multiple household batteries potentially causes legal problems; data protection issues need to be discussed as well. Hence, future research needs to tackle challenges beyond the pure mathematical scope.



\bibliographystyle{IEEEtran}
\bibliography{main}

\appendices

\section{Vector Auto-regressive model}
\label{appendix_VAR}
\subsection{Model Calibration}
There are different methods to calibrate a VAR model \cite{embrechts2011quantitative} \cite{lutkepohl2005new}, here we use maximum likelihood estimation (MLE): Given a data sample the conjoint probability distribution function, called likelihood function here, is optimized with regards to the distribution parameters. As we assume that errors $\epsilon_{t}$ are Gaussian distributed with zero mean this mainly refers to the matrix $A$ in \eqref{eq:VAR} and the corresponding covariance matrix $\Sigma$. The likelihood function $LF(A,\Sigma|\vec{P}^1,\ldots,\vec{P}^n)$  of our VAR(1) model given a data sample $(\vec{P}^1,\ldots,\vec{P}^n)$ reads as follows:
\begin{equation*}
\begin{split}
    LF = & \prod_{t=1}^{n} \frac{1}{(2\pi)^{\frac{k}{2}}\vert \Sigma \vert^{\frac{1}{2}}}\exp\left(-\frac{1}{2}(\epsilon^{\top}_{t}\Sigma^{-1}\epsilon_{t})\right)
\end{split}
\end{equation*}
To simplify optimization, the logarithm is applied to $LF$. The function is then called log-likelihood function, short LLF: 
\begin{equation*}
    \log(LF) = -\frac{nk}{2}\log(2\pi)-\frac{n}{2}\log\left(\vert \Sigma\vert\right) -\frac{1}{2} \sum_{t = 1}^{n}\epsilon^{\top}_{t}\Sigma^{-1}\epsilon_{t}
\end{equation*} 
where $\epsilon_{t} = \vec{P}^t - A \vec{P}^{t-1}$ and $\vert\Sigma\vert$ is the determinant of $\Sigma$. On substituting the value of $\epsilon_t$ in the log-likelihood function and differentiating partially w.r.t to $A,\text{ and } \Sigma$ and setting them to zero an optimum can be found. Note that this method is implemented in all common software applications.
\subsection{Forecast Generation}
The VAR(1) model allows to derive forecasts as well. E.g. for day \( t+1 \), the forecasted price vector \( \hat{\vec{x}}_{t+1} \) is given by:
\begin{equation*}
    \hat{\vec{x}}_{t+1} = \mathbf{A} \vec{x}_t
\end{equation*}
Forecasts for several days ahead can be derived accordingly.
\section{Glossary}
\begin{itemize}
\item \textbf{EPEX Spot}: Spot market platform of the European Energy Exchange. Thereby, spot market refers to short term trading products with physical delivery such as individual hours the next day or 15 minutes blocks the same day.
\item \textbf{Least squares Monte Carlo method}: A valuation method for contracts that include a stochastic price component, particularly for physical assets such as storage's or power plants.
\item \textbf{Maximum likelihood estimation}: A parameter estimation method for a distribution function given a data sample. The likelihood of seeing the given data sample is maximized with regards to the available parameters.
\item \textbf{Time to break even}: The break even point refers to the first moment in an investment's live cycle where the total costs equal the total revenues. The corresponding time measured from the beginning of the investment is called \textit{time to break even}.
\item \textbf{VAR Model}: Vector autoregressive model that allows to capture linear dependencies over time.
\end{itemize}

\end{document}